\documentclass[11pt,twoside]{article}
\usepackage{mathbbold}

\usepackage{amsmath, amssymb, mathrsfs, graphicx}
\usepackage{mathabx}
\usepackage[all]{xy}

\def\az{\alpha}  
    \def\dz{\delta}
    \def\fz{\varphi}
  
\def\lz{\lambda}

\def\Re{\text{\rm Re}}

\def\qd{\quad}
\def\qqd{\qquad}

\newcommand{\mathsym}[1]{{}}

\def\le{\leqslant}
\def\ge{\geqslant}

\font\cms=cmss9 scaled \magstep1

\def\nnd{\noindent}
\def\lmm{\nnd\bg{lmm1}}
\def\prp{\nnd\bg{prp1}}

\def\xmp{\nnd\bg{xmp1}}

\def\algo{\nnd\bg{algo1}}

\def\delmm{\end{lmm1}}
\def\deprp{\end{prp1}}
\def\dexmp{\end{xmp1}}

\def\dealgo{\end{algo1}}

\def\prf{\medskip \noindent {\bf Proof}. }
\def\qed{\text{\quad $\Box$}}
\def\deprf{\qed\medskip}
\def\bg{\begin}
\def\be{\bg{equation}}
\def\de{\end{equation}}

\def\dear{\end{eqnarray}}
\def\lb{\label}
\def\ct{\cite}

\newcommand{\rf}[2]{[\ref{#1}; #2]}

\def\den{\end{enumerate}}

\def\bbb{\mathbbold}

\begin{document}
%\song

%\setcounter{page}{1370}

\thispagestyle{empty}
\renewcommand{\thefootnote}{\fnsymbol{footnote}}

%\vspace*{-0.9in}
%\noindent {Chinese Journal of Applied Probability and Statistics 2015}\newline
%\noindent {January 2013, Volume 29, Issue 1, pp 1-32 }

\begin{center}
{\bf\Large Global algorithms for maximal eigenpair}
\vskip.15in {Mu-Fa Chen}\\
\vskip.15in {Assisted by Yue-Shuang Li}
\end{center}
\begin{center} (Beijing Normal University)\\
\vskip.1in April 29, 2017       %%April 21, 2016
\end{center}
\vskip.1in

\markboth{\sc Mu-Fa Chen}{\sc Maximal Eigenpair}

%\title{Exponential Convergence Rate in Entropy}

%\author{Mu-Fa Chen}

\date{}
%{February 24, 2013}

%\maketitle

%\footnotetext{Received May 17, 2012; accepted June 18, 1012}
%\footnotetext{2000 {\it Mathematics Subject Classifications}.\quad 26D10, 60J60, 34L15.}
%\footnotetext{{\it Key words and phrases}.\quad
%Hardy-type inequality, optimal constant, variational formulas, approximating procedure.}

\begin{abstract} This paper is a continuation of \ct{cmf16} where
an efficient algorithm for computing
the maximal eigenpair was introduced first for tridiagonal matrices and then extended to the irreducible matrices with nonnegative off-diagonal elements. This paper introduces two global algorithms for computing the maximal eigenpair in a rather general setup, including even a class of real (with some negative off-diagonal elements) or complex matrices.
\end{abstract}

\nnd {\small 2000 {\it Mathematics Subject Classification}: 60J60, 34L15}

\nnd {\small {\it Key words and phrases}. Maximal eigenpair, shifted inverse iteration, global algorithm.}

\section{Introduction}\lb{s-01}

To compute the maximal eigenpair of the tridiagonal matrices with positive sub-diagonal elements, an efficient algorithm was introduced \rf{cmf16}{\S  3}.
In the tridiagonal case, the construction of the initials for the algorithm is explicit. In some sense, the results are more or less complete (a modified algorithm, Algorithm \ref{t-17}, is included
in \S 4.4). Next, the algorithm was extended to the general case in \rf{cmf16}{\S  4} which is still efficient for tridiagonally dominant matrices. Note that the initial $v_0$ constructed in \rf{cmf16}{\S  4.2} may not be efficient enough, since the shape of the maximal eigenvector
can be rather arbitrary, could be quite far away from $v_0$ constructed in \rf{cmf16}{\S 4.2}. Thus, we are worrying about the efficiency of the extended algorithm and moreover a global
algorithm is still missed in our general setup. This is the aim of this paper.
In \S  \ref{s-03}, a part of the off-diagonal elements of the matrices are allowed to be negative. We can even handle with some complex matrices.
Let us concentrate on the nonnegative matrices from now on, unless otherwise is stated.

By a shift if necessary, unless otherwise stated, we assume that the given matrix $A=(a_{ij}: 0\le i, j\le N)$
is irreducible and nonnegative: $a_{ij}\ge 0$. We now state our algorithms.
To guarantee the convergence of the iterations in the paper, we assume
that the matrix is irreducible having positive trace, or equivalently,
\be A^n>0 \qqd\text{for each } n\ge \text{ some }n_0.   \lb{01}\de
We mention that in the present nonnegative case, the condition having positive
trace is not serious, otherwise, simply adopt a shift as mentioned at the beginning
of \ct{cmf16}.

In what follows, we omit, without mention time by time, the trivial
case that $\sum_j a_{ij}\!\equiv$\,constant $m\!>\!0$. Since then the maximal
eigenpair of $A$ becomes $(m, \bbb{1})$, where $\bbb{1}$ is the constant
function having components $1$ everywhere.

Recall that the choice of the initials is quite essential for the Rayleigh Quotient Iteration
(RQI), a special shifted inverse iteration. In general, it seems no hope at the
moment to have such explicit analytic
formulas as used in \rf{cmf16}{\S  3}. Instead, as suggested in many
textbooks, one may use other approach to
obtain in a numerical way the required initials, say use the power iteration
for instance. The last approach is safe, but rather slow as shown at the beginning
of \ct{cmf16}. This leads us to come back to the shifted inverse iterations which is a fast cubic algorithm. The ratio of the numbers of iterations for these two algorithms can be thousands.
Now, a critical point is
to avoid the dangerous pitfalls, i.e., the region $(0, \rho(A))$, where $\rho(A)$
is the maximal eigenvalue of $A$. The answer is given in part (1) of the
next two algorithms. At the moment, we are interesting in the generality and
safety, do not take care much about the convergence speed, perhaps, maybe some price
we have to pay here. We will see soon what happen in the next section.

\algo\lb{t-01}{\rm(Specific Rayleigh quotient iteration)}\qd{\cms Let $A=(a_{ij})$ be given.
\bg{itemize} \setlength{\itemsep}{-0.6ex}
\item[(1)] Define column vectors
$$w^{(0)}=(1, 1, \ldots, 1)^*,\qqd  v^{(0)}=w^{(0)}\big/\sqrt{N+1},$$
and set
$$z^{(0)}=\max_{0\le i \le N} \big(A w^{(0)}\big)_i.$$
\item[(2)] For given $v: = v^{(n-1)}$ and $z:=z^{(n-1)}$, let $w:=w^{(n)}$
solve the equation
 \be   (z I-A) w=v. \lb{02}
 \de
As in step $(1)$, define $v^{(n)}=w\big/\sqrt{w^*w}$. Next, define
$$x^{(n)}=\min_{0\le j\le N} \frac{(A w^{(n)})_j}{w_j^{(n)}},\qd
y^{(n)}=\max_{0\le j\le N} \frac{(A w^{(n)})_j}{w_j^{(n)}}, \qd
z^{(n)}= {v^{(n)}}^*A v^{(n)}. $$
\item[(3)] If at some $n\ge 1$, $y^{(n)}-x^{(n)}<10^{-6}$
(or $|z^{(n)}-z^{(n+1)}|<10^{-6}$)(say!), then stop the computation.
At the same time, regard $\big(z^{(n)}, v^{(n)}\big)$ as an approximation of the maximal eigenpair.
\end{itemize}
}\dealgo

The algorithm was presented in \rf{cmf16}{\S  4.1: Choice I}. The simplest choice $v_0$ is reasonable
in the sense that it enables us to cover the general case. We did not pay enough attention on this algorithm
since it looks less efficient. However, as some examples will be illustrated below, this algorithm is actually rather powerful. It is the place to state the main new algorithm of the paper.

\algo\lb{t-02}{\rm(Shifted inverse iteration)}\qd{\cms Everything is the same as in Algorithm \ref{t-01},
except $y^{(n)}$ and $z^{(n)}$ defined in parts $(2)$ and $(3)$ there are exchanged. Moreover, the
resulting $z^{(n)}$ (resp. $x^{(n)}$) is decreasing (resp. increasing) in $n$.
}\dealgo

Let us repeat the sequences $z^{(n)}$, $y^{(n)}$ and $x^{(n)}$ defined in Algorithm \ref{t-02}:
$$x^{(n)}=\min_{0\le j\le N} \frac{(A w^{(n)})_j}{w_j^{(n)}},\qd
y^{(n)}= {v^{(n)}}^*A v^{(n)},\qd
z^{(n)}=\max_{0\le j\le N} \frac{(A w^{(n)})_j}{w_j^{(n)}}. $$
It is obvious that
$$x^{(n)}\le y^{(n)}\le z^{(n)}.$$
In general, Algorithm \ref{t-01} is a little effective than Algorithm \ref{t-02},
saving one iteration for instance, but in Algorithm \ref{t-02}, each iteration is
safe, never failed into the pitfall. This is based on the following dual variational formula.

\prp{\rm \rf{nout08}{Theorem (8)}}\lb{t-11}\;\;{\cms For a nonnegative irreducible matrix $A$, the Collatz--Wielandt formula holds:
$$\sup_{x>0}\min_{i\in E} \frac{(Ax)_i}{x_i}=\rho (A)
= \inf_{x>0}\max_{i\in E} \frac{(Ax)_i}{x_i}.$$
}\deprp

\nnd Actually, suppose that we have $w^{(n-1)}>0$ in Algorithm \ref{t-02}.
Then by Proposition \ref{t-11}
and step (2) of Algorithm \ref{t-02}, we have
$z^{(n-1)}>\rho(A)$ and then the solution $w^{(n)}$ to the equation (\ref{02}) should be
positive: $w^{(n)}>0$. Otherwise, if $z^{(n-1)}<\rho(A)$, then the solution $w^{(n)}$ is negative.
This is the main reason why we choose such a $z^{(n-1)}$ for each $n\ge 1$ in
Algorithm \ref{t-02} and in the case of $n=0$ in Algorithm \ref{t-01} as
our shift, avoiding the change of signs. Note that in Algorithm \ref{t-01}
we adopt $y^{(n)}$ at each step $n\ge 1$, hence the solution
$w^{(n)}$ changes its sign often. This seems dangerous
because $y^{(n)}$ is located in the dangerous region, but up to now,
we have not meet trouble. Therefore, it is still regarded as one
of our two main algorithms. Certainly, if necessary, you can replace
$z^{(k)}$ defined in Algorithm \ref{t-01} by $z^{(k)}=\max_{0\le j\le N} (A w^{(k)})_j/w^{(k)}_j$ for some $k=1, 2, 3$ or so. Or, simply use Algorithm \ref{t-02} instead.

A careful comparison of Algorithm \ref{t-01} and the powerful one
introduced in \rf{cmf16}{\S  3} is delayed to the Appendix.

An easier way to see the efficiency of Algorithms \ref{t-01} and \ref{t-02} is comparing them
with the one given in \rf{cmf16}{\S  4.2}. Suppose that we have used
three iterations in computing a model using the method introduced in
\rf{cmf16}{\S  4.2}, this means on the one hand we have solved the
linear equations in three times. On the other hand, we have solved three
more times in advance to figure out the initials $v^{(0)}$ and $z^{(0)}$
in terms of the triple $(\psi, h, \mu)$. Altogether, we have solved
six linear equations. Or in other words, we have used 6 iterations
in the computation for the specific model. Thus, Algorithms \ref{t-01} and \ref{t-02} should be
regarded as efficient one if no more than 6 iterations are used in the
computation for the same model. As we will see soon, we are actually in
such a successful situation.

To conclude this section, we rewrite Algorithms \ref{t-01} and \ref{t-02} to
a class of matrices with nonnegative off-diagonal elements and negative diagonal
elements: $Q=(q_{ij})$:
$$q_{ij}\ge 0,\;\; i\ne j; \qqd \sum_{j=0}^N q_{ij}\le 0,\;\; 0\le i\le N.$$
In this case, we are studying the maximal eigenpair of $Q$, or alternatively,
the minimal eigenpair of $-Q$. To which, the next two algorithms are devoted.

Again, the trivial case that $\sum_{j=0}^N q_{ij}$ equals a constant is
ignored throughout the paper.

\algo\lb{t-03}{\rm(Specific Rayleigh quotient iteration)}\qd{\cms Let $Q=(q_{ij})$ be given.
\bg{itemize} \setlength{\itemsep}{-0.6ex}
\item[(1)] Define column vectors
$$w^{(0)}=(1, 1, \ldots, 1)^*,\qqd  v^{(0)}=w^{(0)}\big/\sqrt{N+1},$$
and set $z^{(0)}=0.$
\item[(2)] For given $v: = v^{(n-1)}$ and $z:=z^{(n-1)}$, let $w:=w^{(n)}$
solve the equation
 \be   (-Q-z I) w=v. \lb{03}
 \de
As in step $(1)$, define $v^{(n)}=w\big/\sqrt{w^*w}$. Next, define
$$\!x^{(n)}\!\!=\!\min_{0\le j\le N}\! \!\frac{((-Q) w^{(n)})_j}{w_j^{(n)}},\;
y^{(n)}\!\!=\!\max_{0\le j\le N}\!\! \frac{((-Q) w^{(n)})_j}{w_j^{(n)}}, \;
z^{(n)}\!\!=\! {v^{(n)}}^*\!(-Q) v^{(n)}. $$
\item[(3)] If at some $n\ge 1$, $y^{(n)}-x^{(n)}<10^{-6}$ (or $|z^{(n)}-z^{(n+1)}|<10^{-6}$)(say!), then stop the computation.
At the same time, regard $\big(z^{(n)}, v^{(n)}\big)$ as an approximation
of the minimal eigenpair.
\end{itemize}
}\dealgo

\algo\lb{t-04}{\rm(Shifted inverse iteration)}\qd{\cms Everything is the same as in Algorithm \ref{t-03},
except $x^{(n)}$ and $z^{(n)}$ defined in parts $(2)$ and $(3)$ there are exchanged. Moreover, the
resulting $z^{(n)}$ (resp. $x^{(n)}$) is increasing (resp. decreasing) in $n$.
}\dealgo

Algorithms \ref{t-03} and \ref{t-04} are based on \rf{cmf16}{Corollary 12}, a
corollary of Proposition \ref{t-11}.

\section{Examples}

To illustrate the power of the algorithms introduced in the last section,
we examine some typical examples in this section.

To go to practical computation for concrete models, our readers are urged
to prepare enough patience, one may have a large number of iterations since the
initials given in part (1) are quite rough.

The efficient application of Algorithm \ref{t-01} was illustrated by
\rf{cmf16}{Examples 13--16}. To have a concrete comparison of the present algorithms with the one introduced in \rf{cmf16}{\S  4.2}, let us consider a simple example.

\xmp{\rm\rf{cmf16}{Example 21}}\lb{??}\;\;{\rm Let
$$Q=\begin{pmatrix}
-3 & 2 & 0 & 1 & 0\\
4 &  -7 &  3 & 0 & 0\\
0 & 5 & -5 & 0 & 0\\
10 & 0 & 0 & -16 &  6\\
0 & 0 & 0 & 11 & -11-b_4
\end{pmatrix}.  $$
Corresponding to different $b_4$, the minimal eigenvalue $\lz_0$
of $-Q$ and its approximation are as follows. Here and in what follows,
we stop at $z^{(2)}$ once the outputs $z^{(k)}=z^{(2)}$ for every $k\ge 2$.
%\vspace{-0.25truecm}
\begin{center}Table 1. The outputs by Algorithm \ref{t-01}\end{center}
\vspace{-0.2truecm}
$$\begin{tabular}{c|c|c|c}
\hline
$b_4$  & ${\displaystyle z^{(1)}}$ & $z^{(2)}$ &$z^{(3)}=\lz_{\min}(-Q)$ \\
 \hline\hline
0.01  & 0.000278773  &0.000278686 & {  } \\
 \hline
 1  & 0.0251531  &  0.0245175 & {  } \\
 \hline
100  & 0.191729 &  0.182822 &  0.182819 \\
 \hline
 $10^4$  & 0.201695  &  0.195019 &  0.195015 \\
 \hline
\end{tabular}$$
\begin{center}Table 2. The outputs by Algorithm \ref{t-02}\end{center}
\vspace{-0.2truecm}
$$\begin{tabular}{c|c|c|c}
\hline
$b_4$  & ${\displaystyle z^{(1)}}$ & $z^{(2)}$ &$z^{(3)}=\lz_{\min}(-Q)$ \\
 \hline\hline
0.01  & 0.000278637  &0.000278686\,$=\lz_{\min}(-Q)$ & {  } \\
 \hline
 1  & 0.0241546  &  0.0245175 & {  } \\
 \hline
100  & 0.168776 &  0.18275 &  0.182819 \\
 \hline
 $10^4$  & 0.179525  &  0.194932 &  0.195015 \\
 \hline
\end{tabular}$$
While the outputs by the algorithm given in \rf{cmf16}{\S  4.2} are the following.  %\vspace{-0.25truecm}
\begin{center}Table 3. The outputs by the algorithm given in \ct{cmf16}\end{center}
\vspace{-0.2truecm}
$$\begin{tabular}{c|c|c|c}
\hline
$b_4$  & ${\displaystyle z^{(1)}}$ & $z^{(2)}$ &$z^{(3)}=\lz_{\min}(-Q)$ \\
 \hline\hline
0.01  & 0.000278573  &0.000278686 & {  } \\
 \hline
 1  & 0.0236258  &  0.0245174 & 0.0245175 \\
 \hline
100  & 0.200058  &  0.182609 &  0.182819 \\
 \hline
\end{tabular}$$
}
\dexmp

These tables show that the three algorithms are more or less at the same
level of effectiveness. However, the first two are actually more economic
since the last one requires an extra work computing the initial $v_0$.

Comparing \rf{cmf16}{Example 15} with the corrected version of \rf{cmf16}{Example 20} and its improvements
given in \rf{cmf16}{Tables 11, 12} (see the author's homepage), we see that the extended algorithm introduced in \rf{cmf16}{\S  4.2} can be less efficient than
Algorithm \ref{t-01}, it has some limitation for general non-symmetrizable (non-symmetric) matrices.
We call a matrix $A=(a_{ji})$ is symmetrizable, if there exists a positive measure $(\mu_i)$ such that
$$\mu_i a_{ij}=\mu_j a_{ji},\qqd i\ne j.$$
A simple necessary condition for the symmetrizability is
$$a_{ij}> 0\Longleftrightarrow a_{ji}>0,\qqd i\ne j.$$
Refer to \rf{cmf04}{Chapter 7} and references within for the solution
to the symmetrizability problem.

Let us start at a class of non-symmetrizable matrices which are taken
from the so-called single birth $Q$-matrix (cf. \ct{cmf04} and references
within). Define
\be Q\!=\!\begin{pmatrix}
-1 & 1 &0& 0&\cdots\cdots & 0 &0\\
a_1 & -a_1\!-\!2 &2&0& \cdots\cdots & 0 &0\\
a_2 &0\;\; & -a_2\!-\!3&3& \cdots\cdots & 0 &0\\
\vdots& \vdots &\vdots&\vdots &\cdots\cdots &N-2 &0\\
a_{N-1} & 0 &0& 0&\cdots &\!\!\!\! -a_{N-1}\!-\!N\!+\!1 & N\!-\!1\\
a_N & 0 &0& 0&\cdots\cdots & 0 &\!\!-a_N\!-\!N
\end{pmatrix}.  \lb{04-1}\de
For this matrix, we have computed several cases:
$$a_k=1/(k+1),\;\; a_k\equiv 1,\;\; a_k=k,\;\; a_k=k^2.$$
Among them, the first one is hardest and is hence presented below.

\xmp{\rm Let $Q$ be defined by (\ref{04-1}). For different $N$,
the outputs of Algorithm \ref{t-04} (equivalently, Algorithm \ref{t-02})
are as follows.   %% Single_birth2_no_moving.m
\begin{center}Table 4. The outputs for different $N$ by Algorithm \ref{t-04}\end{center}
\vspace{-0.2truecm}
$$\begin{tabular}{c|c|c|c|c|c|c}
\hline
$N$  & ${\displaystyle {z^{(1)}}}$ & $z^{(2)}$ &$z^{(3)}$ &$z^{(4)}$ &$z^{(5)}$ &$z^{(6)}$\\
 \hline\hline
$8$ & $ 0.276727$  & $0.427307$ & $0.451902$ &0.452339 & ${  } $ & ${  } $\\
 \hline
$16$ & $ 0.222132$  & $ 0.367827$ & $0.399959$ & 0.400910& ${  } $ &{  }\\
 \hline
$32$ & $ 0.187826$  & $0.329646$ & $0.370364$ & $0.372308$ & 0.372311 &{  }\\
 \hline
$50$ & $ 0.171657$  & $0.311197$ & $0.357814$ & ${0.360776} $ & 0.360784 &{  }\\
 \hline
$100$ & $0.152106$  & $0.287996$ & $0.343847$ & $0.349166$ & 0.349197 &{  }\\
 \hline
$500$ & $0.121403$  & $0.247450$ & $0.321751$ & $0.336811$ &0.337186 &{ } \\
 \hline
$1000$ & $0.111879$  & $0.233257$ & $0.313274$ & $0.334155 $ &0.335009 & 0.335010 \\
 \hline
$5000$ & $ 0.0947429$  & $0.205212$ & $0.293025$ & $0.328961$ & 0.332609 & 0.332635\\
 \hline
$10^4$ & $0.0888963$  & $0.194859$ & $0.284064$ & $0.326285$ & 0.332113 & 0.332188\\
 \hline
\end{tabular}$$
}
\dexmp
The last line shows that when $N=10^4$, $\lz_{\min}(-Q)\approx 0.332188$. If we
use the shifted matrix $A=Q+mI$, then $\rho(A)\approx 9999.67$. From which, we
get $\lz_{\min}(-Q)\approx 10^4+10^{-4}-9999.67$. Clearly, the second approach
has a less precise output. That is the main difference between Algorithms \ref{t-01}, \ref{t-02}
and \ref{t-03}, \ref{t-04}, even though they are equivalent analytically.

It should be meaningful to have a comparison of the present results with those produced by \rf{cmf16}{\S  4.2}. The outputs listed in the table below come from
the algorithm without using $\dz_1$ defined in that section. For the outputs using $\dz_1$, one more iteration is needed for those $N$ from $16$ to $100$ listed in the table.%%General_Single_birth_2.m
\pagebreak

\begin{center}Table 5. The outputs for different $N$ by the algorithm given in \rf{cmf16}{\S 4.2}\end{center}
\vspace{-0.7truecm}
$$\begin{tabular}{c|c|c|c}
\hline
$N$  & ${\displaystyle {z^{(1)}}}$ & $z^{(2)}$ &$z^{(3)}$ \\
 \hline\hline
$8$ & $ 0.450694$  & $0.452338$  &0.452339 \\
 \hline
$16$ & $ 0.399520$  & 0.400910& ${  } $\\
 \hline
$32$ & $ 0.371433$ & 0.372311 &{  }\\
 \hline
$64$ & $ 0.355722$  & $0.355940$ & ${ }$ \\
 \hline
$100$ & $0.349501 $  & 0.349197 &{  }\\
 \hline
$500$ & $0.340666$  & $0.337185$  &0.337186 \\
 \hline
$1000$ & $0.340871$  & $0.335003$ & 0.335010 \\
 \hline
$5000$ & $ 0.347505$  & $0.332536$  & 0.332635\\
 \hline
$10^4$ & $0.352643$  & $0.331975$ & 0.332188\\
 \hline
\end{tabular}$$
Clearly, the general algorithm introduced in \rf{cmf16}{\S  4.2}
is efficient for this non-symmetrizable model. We have seen that the
present algorithms require more iterations than the earlier one, this
is reasonable since
the computations of the initials are excluded from the last table.
Actually, the computations of the last table cost double time than
the previous one.

The next example is motivated from the classical branching process.
Denote by $(p_k: k\ge 0)$ a given probability measure with $p_1=0$. Let
$$Q\!=\!\left(\!\begin{array}{ccccccc}
-1\!\!\! & p_2 & p_3 &p_4&\cdots\cdots &p_{N-1} &\sum_{k\ge N}p_k \\
2 p_0 &\!\!\! -2 &2 p_2&2p_3 &\cdots\cdots &2 p_{N-2} & 2\sum_{k\ge N-1}p_k \\
0& 3p_0 &\!\!\! -3 & 3p_2&\cdots &3 p_{N-3} & 3\sum_{k\ge N-2}p_k\\
\vdots &\vdots &\vdots &\ddots&\ddots &\ddots &\ddots \\
\vdots &\vdots &\vdots &\ddots&\ddots &-(N-1) & (N-1)p_2\\
0& 0& 0& 0 &\cdots \cdots&\; Np_0\;\; & - Np_0
\end{array}\!\!\right)\!,$$
In the original model, the state $0$ is an absorbing one. Here we regard
it as a killing boundary. Hence it is ruled out from our state space.
Thus, the matrix is defined on $E:=\{1, 2, \ldots, N\}$. Set
$M_1=\sum_{k\in E}k p_k$. When $N=\infty$, in the subcritical case that
$M_1<1$, with a little modification at 0,
it is known that the process generated by $Q$ is ergodic, and is indeed
exponentially ergodic (cf. \rf{rrchen97}{Theorem 1.4\,(iii)}). Hence the exponential convergence rate should
be positive. Otherwise, the
process is not ergodic and so the convergence rate should be zero.

From now on, fix
$$p_0=\az/2,\; p_1=0,\; p_2=(2-\az)/2^2,\; \ldots p_n=(2-\az)/2^n,\cdots,
\qqd \az\in (0, 2).$$
Then $M_1=3(2-\az)/2$ and hence we are in the subcritical case iff $\az\in (4/3, 2)$.

\xmp{\rm Set $\az=1$. Then the outputs of the approximation for
the minimal eigenvalue of $-Q$ by Algorithm \ref{t-02} (or \ref{t-04}) are as follows.   %%revised_GA_4.1.m
\pagebreak

\begin{center}Table 6. The outputs in the supercritical case\end{center}
\vspace{-0.4truecm}
$$\begin{tabular}{c|c|c|c}
\hline
$N$  & ${\displaystyle z^{(1)}}$ & $z^{(2)}$ &$z^{(3)}$ \\
 \hline\hline
$8$ & $ 0.0311491$  & $0.0346044$ & $0.0346310 $\\
 \hline
$16$ & $ 0.00256281$  & $ 0.00260088$ & ${  }$\\
 \hline
\end{tabular}$$
When $N\ge 50 $, $z^{(1)}<10^{-6}$. Hence, $z^{(n)}$ decays quite
quick to zero when $N\to\infty$ (for $n\ge 2$). This is reasonable since
we are now away from the subcritical region.
}\dexmp  %% branching_process_4_28.m
\xmp{\rm Set $\az=7/4$. We are now in the subcritical case and so the
maximal eigenvalue should be positive. We want to know how fast the
local maximal eigenvalue becomes stable (i.e., close enough to the
converge rate at $N=\infty$). Again, we adopt Algorithm \ref{t-02} (or \ref{t-04}).
Up to $N=10^{4}$, the steps of the iterations we need are no more
than $6$. To fasten the convergence, we adopt a convex combination,
as we did several times in \ct{cmf16}. Replace the original
$z^{(0)}=\max_{0\le j\le N} (A w^{(0)})_j$ by
$$z^{(0)}=\xi \min_{0\le j\le N} (A w^{(0)})_j + (1-\xi) (v^{(0)})^* A v^{(0)}.
$$
In view of the practice on $N=8$, we make the choice that $\xi=0.69$.
Then the outputs of the approximation of
the minimal eigenvalue of $-Q$ for different $N$ are as follows.
\vspace{-0.4truecm}
\begin{center}Table 7. The outputs in the subcritical case\end{center}
\vspace{-0.4truecm}
$$\begin{tabular}{c|c|c|c|c}
\hline
$N$  & ${\displaystyle {z^{(1)}}}$ & $z^{(2)}$ &$z^{(3)}$ &$z^{(4)}$ \\
 \hline\hline
$8$ & $ 0.637800$  & $0.638153$ & ${  } $  & ${  } $\\
 \hline
$16$ & $ 0.621430$  & $ 0.625490$ & $0.625539$ & ${  } $\\
 \hline
$50$ & $ 0.609976$  & $ 0.624052$ & $0.624997$ & ${0.625000} $\\
 \hline
$100$ & $0.606948$  & $0.623377$ & $0.624991$ & $0.625000 $\\
 \hline
$500$ & $ 0.604409$  & $ 0.622116$ & $0.624962$ & ${0.625000} $\\
 \hline
$1000$ & $0.604082$  & $0.621688$ & $0.624944$ & $0.625000 $\\
 \hline
$5000$ & $0.603817$  & $0.620838$ & $0.62489$ & $0.625000$\\
 \hline
$10^4$ & $0.603784$  & $0.620511$ & $0.624861$ & $0.625000$\\
 \hline
\end{tabular}$$
From the above table, we see that for $N$ varies from $8$
to $10^4$, in each case, we need at most 4 iterations only.
The computation in each case costs no more than one minute.
Besides, starting from $N=50$, the final outputs are all
the same: $0.625$, which then can be regarded as a very good
approximation of the maximal eigenvalue at infinity $N=\infty$.
}\dexmp

Hopefully, we have already shown the power of our algorithms.

\section{A class of real or complex matrices}\lb{s-03}\lb{S-3}
This section is out of the scope of \ct{cmf16} which depends heavily
on probabilistic idea. Thanks are given to the extended
Perron--Frobenius theory (\ct{nout06}--\ct{nout12})
which makes this section possible.

First, we consider the real case. The special case that
 all off-diagonal elements of $A$ are negative has been treated above,
using $-Q$ instead of $A$ here. Thus, we are now mainly interested
in the case that a part of the off-diagonal elements are negative. Again,
we are concentrated in the study of the maximal eigenpair.

\prp\qd{\cms Let $A$ be a real matrix satisfying (\ref{01}). Then Algorithms
\ref{t-01} and \ref{t-02} are available.
}\deprp

\prf By \rf{nout06}{Theorem 2.2}, condition (\ref{01}) implies that the
matrix $A$ possesses the strong Perron--Frobenius property. Hence it has
the maximal eigenvalue $\rho(A)$ which is simple, positive and corresponds to a positive
eigenvector. Besides, by \rf{nout06}{Theorem 2.6}, the Collatz--Wielandt formula
given in Proposition \ref{t-11} holds. These facts are enough to use Algorithms \ref{t-01} and \ref{t-02}.
\deprf

The next simple observation is helpful.

\lmm\lb{t-11-1}{\cms Condition (\ref{01}) holds iff
$$A^k>0\qqd \text{\cms for } k=n_0, n_0+1,\ldots, 2n_0-1.$$
}
\delmm

\prf Given $n\ge n_0$, write
$$n=r n_0 + s$$
for some integer $r\ge 1$ and $s=0, 1, \ldots, n_0-1$. If $r=1$, then the
conclusion holds by assumption. Otherwise, let $r\ge 2$. Then express
$$n=(r-1) n_0 + (n_0+s).$$
It follows that
$$A^n=\big(A^{n_0}\big)^{r-1} A^{n_0+s}>0$$
as required.
\deprf

We now illustrate our algorithms by a simple example.

\xmp {\rm \rf{nout08}{Example (7)}}\;\;{\rm Let
$$A =
\begin{pmatrix} -1& 8& -1\\
8& 8& 8\\
-1& 8& 8
\end{pmatrix}.$$
Then
$$
{A^2= \begin{pmatrix}
66& 48& 57\\
48& 192& 120\\
57& 120& 129
\end{pmatrix}>0},\qqd
{A^3= \begin{pmatrix}
261& 1368& 774\\
1368& 2880& 2448\\
774& 2448& 1935
\end{pmatrix}>0}.$$
By Lemma \ref{t-11-1}, condition (\ref{01}) holds with $n_0=2$.
The eigenvalues of $A$ are as follows.
$$17.5124, \qd -7.4675,\qd 4.95513.$$
The corresponding maximal eigenvector is
$$(0.486078,\; 1.24981,\; 1)$$
which is positive.

Here are the outputs of our algorithms. Both algorithms
are started at $z^{(0)}=24$.
\begin{center}Table 8. The outputs for a matrix with more negative elements\end{center}
\vspace{-0.2truecm}
$$\begin{tabular}{c|c|c}
\hline
$n$  & $z^{(n)}$: Algorithm \ref{t-01} & $z^{(n)}$: Algorithm \ref{t-02} \\
 \hline\hline
$1$ & $ 17.3772$  & $18.5316$ \\
 \hline
$2$ & $ 17.5124$  & $ 17.5416$ \\
 \hline
 $3$ & ${   }$  & $ 17.5124$ \\
 \hline
\end{tabular}$$

}
\dexmp

Next, we turn to study the complex case. Instead of (\ref{01}), we assume
that
\be \Re\big(A^n\big)>0 \qqd\text{for } n\ge \text{ some }n_0.   \lb{04}\de
Certainly, as usual $\Re(A)$ means the real part of a complex matrix $A$.
This condition is based on \rf{nout12}{Theorems 2.3 and 2.2}, from which we know that $A$ has the maximal, simple, positive eigenvalue. Then we have a weak extension of the Collatz--Wielandt formula as follows.

\prp\lb{13}{\rm \rf{nout12}{Theorems 2.3 and 2.4}}\;\;{\cms Let $A^k\ne 0$ for each $k\ge 1$
and $\Re(A^n)\ge 0$ for every large enough $n$. Then we have for each $x>0$
$$\min_{0\le j\le N}\frac{(\Re(A)x)_j}{x_j}\le \rho(A)\le \max_{0\le j\le N}\frac{(\Re(A)x)_j}{x_j}.$$
}\deprp

Since for the complex conjugate ${\bar x}^*$ of $x$,
the quantity ${\bar x}^* A x$ may still be complex, in view of this, Proposition \ref{13} and the positivity of $\rho(A)$ by (\ref{04}), it seems not reasonable to use ${\bar x}^* A x/({\bar x}^* x)$ as a shift. In
this sense, we do not have a modified version of Algorithm \ref{t-01}.
Fortunately, Algorithm \ref{t-02} is still meaningful.

\algo\lb{t-30}{\rm(Shifted inverse iteration)}\qd{\cms Assume (\ref{04}).
\bg{itemize} \setlength{\itemsep}{-0.6ex}
\item[(1)] Define column vectors
$$w^{(0)}=(1, 1, \ldots, 1)^*,\qqd  v^{(0)}=w^{(0)}\big/\sqrt{N+1},$$
and set
$$z^{(0)}=\max_{0\le i \le N} \big(\Re(A) w^{(0)}\big)_i.$$
\item[(2)] For given $v: = v^{(n-1)}$ and $z:=z^{(n-1)}$, let $w:=w^{(n)}$
solve the equation
 \be   (z I-A) w=v. \lb{05}
 \de
As in step $(1)$, define $v^{(n)}=w\big/\sqrt{{\widebar w}^*w}$. Next, define
$$
z^{(n)}=\max_{0\le j\le N} \frac{(\Re(A) \Re(w^{(n)}))_j}{\Re(w^{(n)})_j}, \qqd
y^{(n)}= ({\bar v}^{(n)})^* A v^{(n)}. $$
\item[(3)] If at some $n\ge 1$, $|y^{(n+1)}-y^{(n)}|<10^{-6}$ (say!), then stop the computation.
At the same time, regard $\big(y^{(n)}, v^{(n)}\big)$ as an approximation
  of the maximal eigenpair.
\end{itemize}
}\dealgo

Note that in Algorithm \ref{t-30}, the sequence $\big\{y^{(n)}\big\}_{n\ge 0}$, but not $\big\{z^{(n)}\big\}_{n\ge 0}$, converges
to $\rho(A)$. To illustrate the use of the algorithm, we consider the following example.

\xmp{\rm \rf{nout12}{Example 2.1}}\;\;{\rm Let
$$A=
\begin{pmatrix}
0.75 - 1.125\, i \; & 0.5882 - 0.1471\, i \; &
  1.0735 + 1.4191\, i\\
   -0.5 -  i & 2.1765 + 0.7059\, i &
  2.1471 - 0.4118\, i\\
   2.75 - 0.125\, i &
  0.5882 - 0.1471\, i & -0.9265 + 0.4191\, i
\end{pmatrix},$$
where the coefficients are all accurate, to four decimal digits.
Then $A$ has eigenvalues
$$3,\qd -2 - i,\qd 1 + i$$
with maximal eigenvector
$$(0.408237,\;\; 0.816507,\;\; 0.408237).$$
The outputs of Algorithm \ref{t-30} are as follows.
\begin{center}Table 9. The outputs for a complex matrix\end{center}
\vspace{-0.2truecm}
$$\begin{tabular}{c|c|c}
\hline
$y^{(1)}$  & $y^{(2)}$ & $y^{(3)}$ \\
 \hline\hline
$3.03949 - 0.0451599\, i$ & $ 3.00471 - 0.0015769\,i$  & $3$ \\
 \hline
\end{tabular}$$
}\dexmp

\section{Appendix}

\subsection{Proof of the last assertion in Algorithm \ref{t-02}}

\prp\lb{t-40}{\cms The sequence
$$ z^{(n)}=\max_{0\le j\le N} \frac{(A w^{(n)})_j}{w_j^{(n)}}
\qd \bigg(\text{\cms resp. } x^{(n)}=\min_{0\le j\le N} \frac{(A w^{(n)})_j}{w_j^{(n)}}\bigg) $$
defined in Algorithm \ref{t-02} is decreasing (resp. increasing ) in $n$.}
\deprp

\prf Let $w>0$ and define
$${\bar\rho}=\max_{0\le j\le N}\frac{(A w)_j}{w_j}.$$
Then $(A w)_j\le {\bar\rho} w_j$ for every $j$. That is,
$$(A_z w)_j\le {\bar\rho}_z w_j\;\;\; \forall j,\qqd A_z:=A/z, \; {\bar\rho}_z={\bar\rho}/z.$$
Since $A_z\ge 0$, it follows that
$$A\sum_{n=0}^\infty A_z^n w \le A\bigg(w +{\bar\rho}_z\sum_{n=0}^\infty A_z^n w\bigg)
\le {\bar\rho} w + \sum_{n=1}^\infty {\bar\rho} A_z^n w
={\bar\rho}\sum_{n=0}^\infty A_z^n w.$$
This means that
$$A (I-A_z)^{-1}w\le {\bar\rho}(I-A_z)^{-1}w$$
since $z> \rho(A)$ by assumption and then $\rho(A_z)<1$. Hence
$$\max_{0\le j\le N}\frac{(A ((I-A_z)^{-1}v))_j}{((I-A_z)^{-1}v)_j}\le {\bar\rho},
\qqd v:=w/\sqrt{w^*w}.$$
Regarding $w=w^{(n-1)}$ and $v=v^{(n-1)}$, this gives us
$$z^{(n)}=\max_{0\le j\le N} \frac{(A w^{(n)})_j}{w_j^{(n)}}
\le \bar\rho
=\max_{0\le j\le N}\frac{(A w^{(n-1)})_j}{w^{(n-1)}_j}
=z^{(n-1)}.$$
Here we have assumed that $z^{(n-1)}>\rho(A)$, otherwise, the computation should
be finished at the step $n-1$.
We have thus proved the assertion on $z^{(n)}$. Dually, we have the assertion on
$x^{(n)}$.
\deprf

\subsection{Proof of the last assertion in Algorithm \ref{t-04}}.

Recall the sequence $\{z^{(n)}\}$ used in Algorithm \ref{t-02} is
given in Proposition \ref{t-40}. Denote by $\{{\tilde z}^{(n)}\}$.
Then, by the relation of $Q$ and $A$ used in Algorithm \ref{t-04}:
$A=Q+ m I$, where $m=\max_i \sum_j a_{ij}$. Hence
$$z^{(0)}I -A =-Q-(m-z^{(0)})I.$$
This means not only ${\tilde z}^{(0)}=0$, but also
$$ w^{(1)}=\big(z^{(0)}I -A\big)^{-1} v^{(0)}=\big(-Q-{\tilde z}^{(0)} I\big)^{-1} v^{(0)}=:{\widetilde w}^{(1)},$$
where ${\tilde w}^{(1)}$ is obtained by the first iteration of Algorithm \ref{t-04}.
Furthermore, similar to the proof of \rf{cmf16}{Corollary 12}, we have
$${\tilde z}^{(1)}=\min_i \frac{(- Q {\widetilde w}^{(1)})_i}{{\widetilde w}^{(1)}_i}
=m-\max_i\frac{(A w^{(1)})_i}{w_i^{(1)}}=m- z^{(1)}.$$
Recursively, we obtain the required assertion.
\qed

\subsection{Comparison of Algorithms \ref{t-01} and \ref{t-03} with the one given in \rf{cmf16}{\S 3}}.\lb{s-4-3}

Since Algorithms \ref{t-01} and \ref{t-03} are equivalent, we need only
to compare Algorithm \ref{t-03} with the one given in \rf{cmf16}{\S 3}. The main difference is their initial $(v^{(0)}, z^{(0)})$.
Clearly, the initial $v^{(0)}$ used in \rf{cmf16}{\S 3}
is finer than the one used in Algorithm \ref{t-03}. Hence, we
need only to compare their $z^{(0)}$.

Next, let $v:=v^{(0)}$ be the initial vector used in
\rf{cmf16}{\S 3}. Denote by $w$ be the solution of
the ordinary inverse iteration (that is the first step of
Algorithm \ref{t-03} or equivalently, Algorithm \ref{t-01}):
$$-Q w =v.$$
Then
\be \frac{(-Q w)_j}{w_j}=\frac{v_j}{((-Q)^{-1}v)_j}=I\!I_j(v)^{-1}.
\lb{07}\de
Here in the last equality of (\ref{07}), we have used the first formula in the proof of \rf{cmf16}{Proposition 23}. Hence
\be\inf_j\frac{(-Q w)_j}{w_j}=\inf_j I\!I_j(v)^{-1}.
\lb{08}\de
The right-hand side of (\ref{08}) is just
$\dz_1^{-1}$ used in \rf{cmf16}{\S 3} as its initial $z^{(0)}$.
The left-hand side of (\ref{08}) should be positive, due to the inverse iteration algorithm, it is certainly bigger than $0$ used as the initial $z^{(0)}$ in Algorithm \ref{t-03}.
In conclusion, both initials used in \rf{cmf16}{\S 3}
are better than those used in Algorithm \ref{t-03}. This completes
the comparison of Algorithm \ref{t-03} and the
one given in \rf{cmf16}{\S 3}.

%\newpage

Naturally, this comparison leads to the next remark.

\subsection{Modification of the algorithm defined in
\rf{cmf16}{\S 3}}\lb{s-4-4}

\nnd{\it Step $1$}. By a shift if necessary, we may assume that we
are given a matrix $Q$ having the form
$$Q\!=\!\left(\!\begin{array}{ccccc}
-(b_0+c_0)\!\!\! & b_0 &0&0 &\cdots \\
a_1 &\!\!\! -(a_1 + b_1+c_1) & b_1 &0 &\cdots \\
0& a_2 &\!\!\! -(a_2 + b_2+c_2) & b_2 &\cdots \\
\vdots &\vdots &\ddots &\ddots &\ddots \\
0& 0& \qquad 0 &\quad a_N^{}\;\; & -(a_N^{}+c_N^{})
\end{array}\!\right)\!,$$
where $a_i> 0,\; b_i> 0,\; c_i\ge 0$ but $c_i\not\equiv 0$.
Note that the maximal eigenvalue of $Q$ is shifted from the
original one but the corresponding eigenvector remain the
same.
\smallskip

\nnd{\it Step $2$}. Following \rf{cmf16}{\S 3}, assume for a moment that some of $c_i\, (i=0, 1, \ldots, N-1)$ is positive.
Then, define
$$\gathered
r_0=1+\frac{c_0}{b_0},\;\; r_n=1+\frac{a_n+c_n}{b_n}-\frac{a_n}{b_n r_{n-1}},\qqd 1\le n<N,\\
h_0=1,\;\; h_n=h_{n-1}r_{n-1}=\prod_{k=0}^{n-1}r_k,\qqd 1\le n \le N,
\endgathered
$$
and additionally,
$$h_{N+1}= c_N h_N+ a_N (h_N-h_{N-1}).$$
We remark that in the special case that
 $$c_0=\cdots=c_{N-1}=0,$$
 by induction, it is easy to check that
 $$r_0=\cdots=r_{N-1}=1$$
 and hence
 $$h_0=\cdots=h_{N}=1.$$
 Furthermore, $h_{N+1}=c_N$.
 Thus, in this special case, we simply ignore
 the sequence $\{h_k\}$ but replace $c_N$ by $b_N$.
 Note that here we use all of the three
sequence $(a_k)$, $(b_k)$ and $(c_k)$ given in $Q$ but no extra thing.
 The role of the sequence $\{h_k\}$ is reducing the
 former case to the last special one and keep the
 same spectrum, in terms of the $H$-transform
${\widetilde Q}$:
\be {\widetilde Q}=\text{\rm Diag}(h_i)^{-1}Q\,\text{\rm Diag}(h_i). \lb{09}\de
The maximal eigenpair $(\rho(Q), g)$ is transformed to
 $\big(\rho\big(\widetilde Q\big)\!\!=\!\rho(Q),  \text{\rm Diag}(h_i)^{-1}g\big)$.

\smallskip

\nnd{\it Step $3$}. In view of Step 2 above, it suffices to consider the
following matrix
\be Q\!=\!\left(\!\begin{array}{ccccc}
-b_0\!\!\! & b_0 &0&0 &\cdots \\
a_1 &\!\!\! -(a_1 + b_1) & b_1 &0 &\cdots \\
0& a_2 &\!\!\! -(a_2 + b_2) & b_2 &\cdots \\
\vdots &\vdots &\ddots &\ddots &\ddots \\
0& 0& \qquad 0 &\quad a_N^{}\;\; & -(a_N^{}+b_N^{})
\end{array}\!\right)\!,  \lb{10}\de
where $a_i,\; b_i> 0$. This step is changed from the original, where
everything we are working here is transfer into the original matrix $Q$
rather than the simpler one here. It seems a direct treatment of the
present matrix $Q$ is slightly simpler.

Define the sequence $(\mu_i)$ as usual:
$$\mu_0=1, \;\; \mu_n=\mu_{n-1}\frac{b_{n-1}}{a_n}=\frac{b_0 b_1 \cdots b_{n-1}}{a_1 a_2 \cdots a_n},
\qquad 1\le n\le N.$$
Next, define
\be\fz_n=\sum_{k= n}^N \frac{1}{\mu_k b_k}, \qqd 0\le n\le N.\de
and
\be\dz_1=\max_{0\le n\le N} \bigg[\sqrt{\fz_n}\sum_{k=0}^n \mu_k \sqrt{\fz_k}+ \frac{1}{\sqrt{\fz_n}}\sum_{n+1\le j \le N}\mu_j \fz_j^{3/2}\bigg].  \de

Having these preparations at hand, we can now start our iterations.
\smallskip

\nnd{\it Step $4$}. As in \rf{cmf16}{\S 3}, choose
\be  w^{(0)}=\sqrt{\fz},\qd v^{(0)}=w^{(0)}/\|w^{(0)}\|_{\mu, 2},
\qd z^{(0)}=\dz_1^{-1},\de
where $\|\cdot\|_{\mu, 2}$ denotes the $L^2(\mu)$-norm. Note that
here in the non-symmetric case, the use of the measure $(\mu_i)$ cannot be ignored since in this case, we are based on, $\dz_k$ for instance, the $L^2(\mu)$ setup.
\smallskip

\nnd{\it Step $5$}. For given $v=v^{(n-1)}$ and $z=z^{(n-1)}$, let
 $w=w^{(n)}$ solve the linear equation
\be (-Q-z I) w=v \de
and then define $v^{(n)}=w/\|w\|_{\mu, 2}$. An explicit solution of this $w$
is now available, refer to \rf{cmf17}{Algorithm 3}.
\smallskip

\nnd{\it Step $6$}. At the $k$th ($k\ge 1$) iteration, in addition to the one $(v^{(k)}, -Q v^{(k)})_{\mu, 2}$ used
in \rf{cmf16}{\S 3}, one may also adopt $z^{(k)}=\dz_k^{-1}$:
\be\dz_k=\max_{0\le i\le N}\frac{1}{v^{(k)}_i}
\bigg[\fz_i \sum_{j=0}^i \mu_j v^{(k)}_j+\sum_{i+1\le j\le N}\mu_j \fz_j v^{(k)}_j\bigg].  \lb{15}\de
This is the main new point in the modified algorithm.
Since \rf{cmf10}{Theorems 2.4\,(3), 3.2\,(1) and (3.6)}, we have
$$\dz_k^{-1} \le \lz_{\min}(-Q)\le (v^{(k)}, -Q v^{(k)})_{\mu, 2}\text{ for each $k$ and $n$}.$$
By \rf{cmf16}{Proposition 23} and \rf{cmf10}{Theorem 3.2\,(1)}, we have known that
the sequence $\{\dz_k^{-1}\}$, deduced in the theorem just cited using the approximating eigenvectors obtained by the ordinary inverse
iteration (without shift), is increasing to $\lz_{\min}(-Q)$. It should be clear that the present
sequence $\{\dz_k^{-1}\}$ produced by the advanced shifted inverse iteration should
converge to $\lz_{\min}(-Q)$ more faster.
Thus the new $z^{(k)}\,(k\ge 1)$ not only avoids the dangerous
region but may also accelerate the convergence of the algorithm. Certainly,
the computation of $\dz_k$ needs more work than the one of
$(v^{(k)}, -Q v^{(k)})_{\mu, 2}$.

The use of the quantity (\ref{15}) is motivated from the remark above on ``Comparison of Algorithms \ref{t-01} and \ref{t-03} with the one given in \rf{cmf16}{\S 3}''. The formula (\ref{15}) is a corollary of \rf{cmf10}{Theorem 2.4\,(3)} which depends on the form (\ref{10}) of $Q$. For general $Q$ such the one in Step 1, we do not have an analog
of \rf{cmf10}{Theorem 2.4\,(3)}, and so (\ref{15}) is not applicable in
such a general situation.
\smallskip

\nnd{\it Step $7$}. To go back to the original matrix $A$, denote its
maximal eigenpair by $(\rho(A), g)$. Recall that the matrix $Q$ at the
beginning is obtained from $A$ by a shift: $Q=A-mI$, $m:=\max_i\sum_j a_{ij}$. Let $(z, v)$ be the output from the last iteration in Step 6.
Then we have
\be\rho(A)\approx m-z, \qqd g\approx \text{Diag}(h_i)v.\lb{16}\de

We now summery the above discussions as a modified algorithm.

\algo\lb{t-17}\;\; {\cms For tridiagonal matrix, the Step 1--Step 7 above consist a modified
algorithm of the one introduced in \rf{cmf16}{\S  3}.
}
\dealgo

We are now ready to study a randomly chosen example, introduced to the author by
Tao Tang, to justify the power of our algorithms and also to compare their efficiency.

\xmp\lb{t-18}\;\;{\rm Let
$$A =\begin{pmatrix}
2.334 & 0.9962 & 0 & 0 & 0 & 0\\
0.5142 & 2.6725 & 0.1111 & 0 & 0 & 0\\
0 & 0.2115 & 2.263 & 0.1405 & 0 & 0\\
0 & 0 & 0.8442 & 2.8457 & 0.7595 & 0\\
0 & 0 & 0 & 0.2347 & 2.2257 & 0.0781\\
0 & 0 & 0 & 0 & 0.9837 & 2.1582\end{pmatrix}.$$
Then the eigenvalues of $A$ are
$$3.26753,\; 3.16247,\; 2.40182,\; 2.12632,\; 1.80416,\; 1.73679.$$
The outputs of our algorithms are given in the table below.
\begin{center}Table 10. Comparison of four algorithms\end{center}
\vspace{-0.2truecm}
$$\begin{tabular}{c|c|c|c|c|c}
\hline
Algorithm & ${\displaystyle {z^{(1)}}}$ & $z^{(2)}$ &$z^{(3)}$ &$z^{(4)}$ &$z^{(5)}$ \\
 \hline\hline
Algorithm 1 & $ 3.30193$  & $3.26737$ & $3.26754 $ & $3.26753$ & ${  } $\\
 \hline
Algorithm 2 & $ 3.64033$  & $ 3.32623$ & $3.26937$ & $3.26756$ & $3.26753$\\
 \hline
Algorithm 17a & $ 3.2618$  & $ 3.26752$ & $3.26753$ & ${ } $ & ${ } $\\
 \hline
Algorithm 17b & $3.27947$  & $3.2685$ & $3.26754$ & $3.26753$\\
 \hline
\end{tabular}$$
where the algorithms in the last two lines mean that

Algorithm 17a: take $z^{(k)}=(v^{(k)}, -Q v^{(k)})_{\mu, 2}$ for each $k\ge 1$.

Algorithm 17b: take $z^{(k)}=\dz_k^{-1}$ defined by (\ref{15}) for each $k\ge 1$.
}\dexmp

\prf To apply Algorithm \ref{t-17}, take $m=4.4494$. Then $Q=A-m I$:
$$Q =
\begin{pmatrix}
-2.1154 & 0.9962 & 0 & 0 & 0 & 0 \\
0.5142 & -1.7769 & 0.1111 & 0 & 0 & 0\\
0 & 0.2115 & -2.1864 & 0.1405 & 0 & 0 \\
0 & 0 & 0.8442 & -1.6037 & 0.7595 & 0 \\
0 & 0 & 0 & 0.2347 & -2.2237 & 0.0781\\
0 & 0 & 0 & 0 & 0.9837 & -2.2912
\end{pmatrix}.$$
We have $h=(2.12347,\; 29.3339,\; 453.284,\; 924.514,\; 24961)$.
The $H$-transform of $Q$ becomes
$${\widetilde Q}=
\begin{pmatrix}
-2.1154 & 2.1154 & 0 & 0 & 0 & 0\\
0.242151 & -1.7769 & 1.53475 & 0 & 0 & 0\\
0 & 0.0153104 & -2.1864 & 2.17109 & 0  & 0\\
0 & 0 & 0.0546316 & -1.6037 & 1.54907 & 0\\
0 & 0 & 0 & 0.115072 &  -2.2237 & 2.10863 \\
0 & 0 & 0 & 0 & 0.0364346 & -2.2912
\end{pmatrix}.$$
Then we are ready to use Algorithm \ref{t-17} for the maximal eigenpair of ${\widetilde Q}$ and finally return to the one for $A$ by (\ref{16}).
\deprf

To explain the word ``modified'' in detail, we transfer Algorithm \ref{t-17}
to the one presented in \rf{cmf16}{\S 3}. To do so, we keep the notation $Q$, $\mu$,
$\fz$, $\dz_1$ and so on used in \rf{cmf16}{\S 3}, but add superscript~$\widetilde{ }$ to those notation used in Steps 3, 4 above. Let $\tilde\mu=h^2\mu$ (i.e., $\tilde\mu_i=h_i^2\mu_i$). Then, as mentioned in \rf{cmf16}{\S 5}, the mapping
$f\to{\tilde f}:=f/h$ gives us not only an isometry from $L^2(\mu)$ to $L^2(\tilde\mu)$
\big(i.e., $\|f\|_{\mu, 2}=\|\tilde f\|_{\tilde \mu, 2}$\big), and then also an isospectrum of $Q$ on $L^2(\mu)$ and $\widetilde Q$ on $L^2(\tilde \mu)$:
$$(f, Q f)_{\mu}=(\tilde f, {\widetilde Q}\tilde f)_{\tilde \mu}, \qqd \|f\|_{\mu, 2}=1.$$
Now, from $L^2(\tilde\mu)$ to $L^2(\mu)$, we have
$${\tilde\fz}_n=\sum_{k= n}^N \frac{1}{{\tilde\mu}_k {\tilde b}_k}
\to \sum_{k= n}^N \frac{1}{h_k h_{k+1}\mu_k b_k}=\fz_n, \qqd 0\le n\le N.$$
Here the transform ${\tilde\mu}_k {\tilde b}_k\to h_k h_{k+1}\mu_k b_k$ for each $k\le N-1$ is regular, except the last term in the sum $\big({\tilde\mu}_N {\tilde b}_N\big)^{-1}$, where ${\tilde b}_N$ is actually the element ${\tilde c}_N$ which is
obtained from the transform $Q\to {\widetilde Q}$, and $h_{N+1}$ and $b_N$ are specified in \rf{cmf16}{\S 3} to make the unified expression in the second sum.
We mention here that $h_{N+1}$ is the original paper \ct{cmf16} should be replaced by
$$h_{N+1}= c_N h_N+ a_N (h_N-h_{N-1})$$
since the sequence $(c_i)$ used in \ct{cmf16} and \ct{chzhx} have different sign. Next,
$$\aligned
{\tilde\dz}_1&=\max_{0\le n\le N} \bigg[\sqrt{{\tilde\fz}_n}\sum_{k=0}^n {\tilde\mu}_k \sqrt{{\tilde\fz}_k}+ \frac{1}{\sqrt{{\tilde\fz}_n}}\sum_{n+1\le j \le N}{\tilde\mu}_j {\tilde\fz}_j^{3/2}\bigg]\\
\to\dz_1&=\max_{0\le n\le N} \bigg[\sqrt{\fz_n}\sum_{k=0}^n \mu_k h_k^2 \sqrt{\fz_k}+ \frac{1}{\sqrt{\fz_n}}\sum_{n+1\le j \le N}\mu_j h_j^2 \fz_j^{3/2}\bigg].
\endaligned $$
At the same time,
$$\aligned
&\big(-{\widetilde Q} -{\tilde z}I\big){\tilde w}={\tilde v}\\
&\Longleftrightarrow \big(-\text{\rm Diag}(h)^{-1} Q\, \text{\rm Diag}(h) -{\tilde z}I\big){\tilde w}={\tilde v}\\
&\Longleftrightarrow (- Q  -{\tilde z}I)\text{\rm Diag}(h){\tilde w}=\text{\rm Diag}(h){\tilde v}\\
&\Longleftrightarrow (- Q  -{z}I) w = v.
\endaligned
$$
Here in the last line, $\tilde z$ is replaced by $z$, this is due to the isospectrum:
an lower bound of the spectrum of $-\widetilde Q$ is also the one of $-Q$.
The fact that $\text{\rm Diag}(h){\tilde w}=w$ comes from the definition of our mapping $f\to \tilde f$. Finally, since the isometry, we have $\|w\|_{\mu, 2}=\|\tilde w\|_{\tilde\mu, 2}$. We have thus deduced the algorithm presented in \rf{cmf16}{\S 3} from the modified
one.

\subsection{Modification of the algorithm introduced in \rf{cmf16}{\S 4.2}}

In parallel to \S \ref{s-4-4}, we may introduce a modification of
the algorithm presented in \rf{cmf16}{\S 4.2}. The main idea is:
once we obtain the function $h$, it can be ignored since we can use
the general transform
$\widetilde Q$ defined in (\ref{09}) instead of the original $Q$ to continue
the procedure of the algorithm constructed in \rf{cmf16}{\S 4.2}.
Since this modification is only a mimic of the one for tridiagonal
matrix (\S \ref{s-4-4}), something may be lost. For instance, the
sequence $\{\dz_k^{-1}\}$ formally defined by (\ref{15}) may no longer be
the lower bound of $\lz_{\min}(-Q)$, one has to take care in practice.
\medskip

To conclude this paper, we remark some possible extension of the algorithms given
here to a more general setup. For a larger class of Markov generators, the algorithms are
meaningful. Actually, the Perron--Frobenius property as well as the
the Collatz--Wielandt formula have been generalized by a number of
authors. In particular, the part of the Collatz--Wielandt formula used
in Algorithm \ref{t-04} as $z^{(n)}$ was extended by
\rf{dv75}{$\psi_2(V)$ in the Theorem}. See also \rf{sh84}{(1.1)\,i) and \S 2} and more recently, \rf{abk16}{Theorem 2.1}. Note the difference: we are working on $\lz_{\min}(-L)$ here rather than $\lz_{\max}(L)$ in the cited papers.

In the nonlinear case, the shifted inverse iteration (Algorithms
\ref{t-02} or \ref{t-04}) is more essential, actually Algorithm \ref{t-01} may no longer be applicable since equation (\ref{02}) often has no real
solution. This point is illustrated in \ct{cmf17} where the shift is
based on a generalization of (\ref{15}). In view of \rf{ckc14}{Theorem
2.3 and Corollary 2.5}, it seems that Algorithm \ref{t-02} and  its variations could be applied to a more general setup.
\medskip

\nnd{\bf Acknowledgments}\qd {\small The author thanks Ms Yue-Shuang Li
for her assistance in computing the large matrices using MatLab, and
also pointed out the error on $h_{N+1}$ mentioned in \S 4.4. The author also acknowledges Mr Xu Zhu for constructing Example \ref{t-18} which leads us to find out the error just mentioned.
Research supported in part by
         National Natural Science Foundation of China (No. 11626245),
         the project from the Ministry of Education in China,
and the Project Funded by the Priority Academic Program Development of Jiangsu Higher Education Institutions.
}

\vspace{-0.25truecm}

\nnd {\small
Mu-Fa Chen\\
School of Mathematical Sciences, Beijing Normal University,
Laboratory of Mathematics and Complex Systems (Beijing Normal University),
Ministry of Education, Beijing 100875,
    The People's Republic of China.\newline E-mail: mfchen@bnu.edu.cn\newline Home page:
    http://math0.bnu.edu.cn/\~{}chenmf/main$\_$eng.htm
}

\end{document}